%% file: 2001-11.tex
\documentclass{gtart}

\input gtoutput

\lognumber{165}
\volumenumber{5}\papernumber{11}\volumeyear{2001}
\pagenumbers{335}{345}
\accepted{2 April 2001}
\received{20 February 2001}
\revised{4 March 2001}
\published{2 April 2001}
\proposed{Haynes Miller}
\seconded{Ralph Cohen, Gunnar Carlsson}

\usepackage{amscd,amssymb,amsmath}

\newcommand{\al}        {\alpha}
\newcommand{\bt}        {\beta} 
\newcommand{\Rhe}       {R^\wedge_e}
\newcommand{\zt}        {\zeta}
\newcommand{\Rhat}      {\widehat{R}}
\newcommand{\dl}        {\delta}
\newcommand{\ep}        {\epsilon}
\newcommand{\tht}       {\theta}
\newcommand{\ov}[1]     {\overline{#1}}
\newcommand{\psb}[1]    {[\![#1]\!]}
\newcommand{\xra}       {\xrightarrow}
\newcommand{\Rb}        {\overline{R}}
\newcommand{\colim}  {\operatornamewithlimits{\underset{\longrightarrow}{lim}}}
\newcommand{\st}        {\;|\;}
\newcommand{\sg}        {\sigma}
\newcommand{\ann}      {\operatorname{ann}}
\newcommand{\invlim} {\operatornamewithlimits{\underset{\longleftarrow}{lim}}}
\newcommand{\sse}       {\subseteq}

\renewcommand{\:}{\colon\thinspace}

\newtheorem{theorem}{Theorem}

\newtheorem{lemma}[theorem]{Lemma}
\newtheorem{proposition}[theorem]{Proposition}
\newtheorem{corollary}[theorem]{Corollary}
\theoremstyle{definition}
\newtheorem{remark}[theorem]{Remark}

\newenvironment{Relax}{\relax}{\relax}

\begin{document}
\title{Complex cobordism of involutions}
\author{N\thinspace P Strickland}
\coverauthors{N\noexpand\thinspace P Strickland}
\asciiauthors{N P Strickland}

\address{Department of Mathematics,
University of Sheffield\\Western Bank, Sheffield, S10 2TN,UK}

\email{N.P.Strickland@sheffield.ac.uk}

\begin{abstract}
 We give a simple and explicit presentation of the
 $\mathbb{Z}/2$--equivariant complex cobordism ring, as an algebra over
 the nonequivariant complex cobordism ring.
\end{abstract}

\asciiabstract{
We give a simple and explicit presentation of the
Z/2-equivariant complex cobordism ring.
}

\keywords{Equivariant, complex cobordism, involution}

\primaryclass{55N22}

\secondaryclass{55N91}

\begin{Relax}
\end{Relax}

\maketitle 

\section{Introduction}

Let $A$ be the group of order two, and let $MU_A^*$ be the homotopical
$A$--equivariant complex cobordism ring, as defined
in~\cite{lemast:esh} (for example).  In this note we give a simple and
explicit presentation of $MU_A^*$ as an algebra over the
nonequivariant complex cobordism ring $MU^*$ (which is of course well
understood).

Our construction and proofs are short and elementary so we have
decided to publish them as they stand.  Elsewhere, we will describe
the conceptual background (formal multicurves and equivariant formal
groups) and generalise our results to other compact Lie groups.  In
the general abelian case, we hope to be almost as explicit as in the
present work.  We will also describe $MU_A^*X$ for many
naturally-occurring $A$--spaces $X$, using the algebraic geometry of
formal schemes; from the right viewpoint, the results are closely
parallel to their nonequivariant analogs.  A large part of the theory
depends only on the fact that $MU_G^*PV$ is a finitely generated
projective module over $MU_G^*$ for all representations $V$, which
holds when $G$ is abelian; we do not know whether there are any
nonabelian examples.

Our approach is based on a pullback square exhibited by
Kriz~\cite{kr:zec}.  We note that Sinha~\cite{si:cce} has also
calculated $MU_A^*$, but his presentation is less explicit and depends
on some choices.

\section{The model}

Let $L=MU^*$ be the Lazard ring, with universal formal group law
$x+_Fy=\sum_{ij}a_{ij}x^iy^j$.  Recall that
\begin{align*}
 a_{ij} &= a_{ji} \\
 a_{0i} &= \dl_{1i} = \begin{cases} 1 & \text{ if } i=1 \\
                                    0 & \text{ otherwise. }
                      \end{cases}
\end{align*}
There are of course other relations, expressing the associativity of
$F$.  It is well-known that $L$ can be expressed as a polynomial
algebra $\mathbb{Z}[x_1,x_2,\ldots]$, although there is no system of
generators that is both explicit and convenient.  Here we take the
structure of $L$ as given and concentrate on describing $MU_A^*$ as an
algebra over it.

Let $R$ be generated over $L$ by elements $s_{ij}$ (for $i,j\geq 0$)
and $t_i$ (for $i\geq 0$) subject to relations given below.  We use
$b_k$ as a synonym for $s_{0k}$, and $e$ as a synonym for
$b_0=s_{00}$.  The relations are:
\begin{align*}
 t_0             &= 0              \\
 s_{10}          &= 1              \\
 s_{i0}          &= 0 \qquad \text{ for } i>1 \\
 t_k    - b_k    &= e t_{k+1} \\
 s_{jk} - a_{jk} &= e s_{j+1,k}
\end{align*}
We can give $R$ a grading with $|s_{ij}|=|a_{ij}|=2(1-i-j)$ and
$|t_k|=2(1-k)$ so $|b_k|=2(1-k)$ and $|e|=2$.

Note that the equation $t_0-b_0=et_1$ gives $(1+t_1)e=0$.  Moreover,
we have $t_1=b_1+et_2=a_{01}+es_{11}+et_2=1+e(s_{11}+t_2)$, so
$1-t_1=0\pmod{e}$.  It follows that $1-t_1^2=(1-t_1)(1+t_1)=0$, so
$t_1^2=1$.

Our main result is as follows:
\begin{theorem}\label{thm-main}
 There is an isomorphism $R\simeq MU_A^*$ of graded $MU^*$--algebras.
\end{theorem}
The rest of this paper constitutes the proof.

\begin{remark}
 Greenlees has studied the ring $L_A$ that classifies $A$--equivariant
 formal group laws; he shows that there is a surjective map
 $\al\:L_A\xra{}MU_A^*$ whose kernel is $e$--divisible and $e$--torsion
 (and is conjectured to be zero).  One can deduce from the above
 theorem that there is a ring map $\bt\:R=MU_A^*\xra{}L_A$ with
 $\al\bt=1$.  It seems likely that our methods can be used to prove
 that $\al$ is an isomorphism, but we have not succeeded as yet.
\end{remark}

\section{The pullback square}

Put $R'=L[b_0,b_1,b_2,\ldots][b_0^{-1}]$.  Next, we define as usual
$[2](x)=x+_Fx=\sum_{i,j}a_{ij}x^{i+j}\in L\psb{x}$, and we put
$\Rhat=L\psb{e}/[2](e)$, and $(\Rhat)'=\Rhat[1/e]$.  There is an
evident map $\Rhat\xra{}(\Rhat)'$, which we call $\zt$.  Next, we
define a map $\xi\:R'\xra{}(\Rhat)'$ by letting $\xi(b_i)$ be the
coefficient of $x^i$ in the power series $x+_Fe$.  (In particular,
$\xi(b_0)=e$, which is invertible in $(\Rhat)'$, as required.)

The following result summarises much of what was previously known
about $MU_A^*$:
\begin{theorem}\label{thm-pullback}
 There is a pullback square as follows:
$$ \begin{CD}
MU_A^* @>>> R' \\
@VVV @VV\xi V\\
\Rhat @>\zt>>(\Rhat)' 
% \node{MU_A^*} \arrow{e} \arrow{s} \node{R'} \arrow{s,r}{\xi} \\
%  \node{\Rhat} \arrow{e,b}{\zt} \node{(\Rhat)'}
 \end{CD}$$
\end{theorem}
\begin{proof}
 This is Theorem~1.1 of~\cite{kr:zec}, in the case $p=2$.  Our element
 $b_i$ is Kriz's $b_1^{(i)}u_1$, and his $u_1$ is our $e$.  There is a
 slight misprint in the statement: the element $b_k^{(0)}$ is $1$, not
 $0$.  Geometrically, the square comes from the Tate diagram of
 cofibrations, as studied in~\cite{grma:gtc}.  In this context,
 $\Rhat$ is identified with $MU^*BA$, and $R'$ with
 $\pi_{-*}\Phi^AMU_A$. 
\end{proof}

Below we shall exhibit a similar pullback square with $MU_A^*$
replaced by $R$; Theorem~\ref{thm-main} follows by uniqueness of
pullbacks.  To show that our square is a pullback, we use the
following result.

\begin{theorem}\label{thm-hasse}
 Let $S$ be a ring, and $d$ an element of $S$.  Suppose that $S$ has
 bounded $d$--torsion, or in other words that
 $\bigcup_k\ann_S(d^k)=\ann_S(d^N)$ for some $N$.  Then there is a
 pullback square as follows:
 $$\begin{CD}
S @>>> S[1/d]\\
@VVV @VVV\\
S^\wedge_d @>>> (S^\wedge_d)[1/d].
%  \node{S} \arrow{e} \arrow{s} \node{S[1/d]} \arrow{s} \\
%  \node{S^\wedge_d} \arrow{e} \node{(S^\wedge_d)[1/d]}
 \end{CD}$$
\end{theorem}
\begin{proof}
 This is standard, but we give a proof for completeness.  First, it is
 clear that there is a commutative square as shown.  Thus, if we let
 $T$ denote the pullback of $S[1/d]$ and $S^\wedge_d$, we get a map
 $S\xra{}T$, and we must show that it is an isomorphism.

 Suppose that $u\in S$ has image zero in both $S[1/d]$ and in
 $S^\wedge_d$.  As $u\mapsto 0$ in $S[1/d]$ we have $d^ku=0$ for some
 $k\geq 0$.  As $u\mapsto 0$ in $S^\wedge_d=\invlim_jS/d^j$ we see
 that $u\mapsto 0$ in $S/d^N$, so $u=d^Nv$ say.  It follows that
 $d^{N+k}v=0$ but $\ann_S(d^{N+k})=\ann_S(d^N)$ so $d^Nv=0$ or in
 other words $u=0$.  This means that the map $S\xra{}T$ is injective.

 Now suppose we have an element $(u,v)\in T$, so $u\in S[1/d]$ and
 $v\in S^\wedge_d$, and $u$ and $v$ have the same image in
 $(S^\wedge_d)[1/d]$.  We then have $u=u'/d^i$ for some $u'\in S$ and
 $i\geq 0$, and $d^ju'=d^{i+j}v$ in $S^\wedge_d$ for some $j\geq 0$.
 Next, we can choose a sequence $(v_m)$ in $S$ such that
 $v_m=v_{m+1}\pmod{d^m}$ and $v_m\xra{}v$ in $S^\wedge_d$.  It follows
 that the element $d^j(u'-d^iv_m)=d^{i+j}(v-v_m)$ maps to zero in
 $S/d^{i+j+m}$, say $d^j(u'-d^iv_m)=d^{i+j+m}w_m$.  We may replace
 $v_m$ by $v_m+d^mw_m$ and thus assume that $d^ju'=d^{i+j}v_m$ for all
 $m$.  Now let $x_m$ be such that $v_{m+1}-v_m=d^mx_m$.  As
 $d^{i+j}v_m=d^ju'=d^{i+j}v_{m+1}$, we see that $d^{i+j+m}x_m=0$.  As
 $\ann_S(d^k)\sse\ann_S(d^N)$ for all $k$, we have $d^Nx_m=0$.  When
 $m\geq N$ this gives $d^mx_m=0$ and thus $v_m=v_{m+1}$.  We put
 $v'=v_N=v_{N+1}=\ldots\in S$.  It is clear that $v'\mapsto v$ in
 $S^\wedge_d$.  Moreover, the equation $d^ju'=d^{i+j}v_m$ shows that
 $v'\mapsto u'/d^i=u$ in $S[1/d]$, so $v'\mapsto (u,v)\in T$.  Thus,
 the map $S\xra{}T$ is also surjective, as required.
\end{proof}
It will be a significant part of our task to check that $R$ has
bounded $e$--torsion.

\section{Completion at $e$}

The following proposition can be proved easily by reading the relations
modulo $e$.  Here $\dl_{1k}$ is Kronecker's delta.
\begin{proposition}
 There is a map $\ep\:R\xra{}L$ with
 \begin{align*}
  \ep(s_{ij}) &= a_{ij}   \\
  \ep(t_k)    &= \ep(b_k) = \dl_{1k} \\
  \ep(e)      &= 0,
 \end{align*}
 and this induces an isomorphism $R/e\xra{}L$. \qed
\end{proposition}

This can be improved as follows:
\begin{proposition}\label{prop-phi}
 There is a map $\phi\:R\xra{}\Rhat$ given by
 \begin{align*}
  \phi(s_{ij}) &= \sum_{k\geq 0} a_{i+k,j}e^k   \\
  \phi(t_k)    &= \sum_{i,j\geq 0} a_{i+k,j}e^{i+j} \\
  \phi(b_l)    &= \sum_{m\geq 0} a_{lm}e^m   \\
  \phi(e)      &= e.
 \end{align*}
 If $\ep'\:\Rhat\xra{}L$ is given by $\ep'(e)=0$, then
 $\ep'\phi=\ep$.  Moreover, $\phi$ induces an isomorphism
 $\hat{\phi}\:\Rhe\xra{}\Rhat$.
\end{proposition}
\begin{proof}
 First, note that the definitions of $\phi(b_l)$ and $\phi(e)$ are
 simply the appropriate specialisations of the definition of
 $\phi(s_{ij})$ so they can essentially be ignored.

 To see that $\phi$ respects the relations, note that
 \begin{align*}
  \phi(t_0)    &= \sum_{ij}a_{ij}e^{i+j}=[2](e)=0 \\
  \phi(s_{10}) &= \sum_{k\geq 0} a_{1+k,0}e^k
                = \sum_{k\geq 0} \dl_{k0}e^k = 1 \\
  \phi(t_k) - \phi(b_k) &= 
   \sum_{i>0}\sum_{j\geq 0} a_{i+k,j}e^{i+j} \\
   &= e \sum_{i,j\geq 0} a_{i+k+1,j}e^{i+j} = e \phi(t_{k+1}) \\
  \phi(s_{ij}) - a_{ij} &= 
   \sum_{k>0} a_{i+k,j}e^k \\
   &= e \sum_{k\geq 0} a_{i+k+1,j}e^k = e \phi(s_{i+1,j}).
 \end{align*}
 It is trivial to check that $\ep'\phi=\ep$.  

 There is an evident map $\psi\:L[e]\xra{}R$, which induces a map
 $\hat{\psi}\:L\psb{e}\xra{}\Rhe$.  As $\Rhat$ is already complete at
 $e$, we also get a map $\hat{\phi}\:\Rhe\xra{}\Rhat$.  It is easy to
 see that the composite
 $\hat{\phi}\hat{\psi}\:L\psb{e}\xra{}\Rhat=L\psb{e}/[2](e)$ is just
 the usual quotient map.

 Next, one checks by induction on $m$ that in $R$ we have
 \begin{align*}
  s_{ij} - \sum_{l<m}a_{i+l,j}e^l &= s_{i+m,j}e^m \\
  t_k    - \sum_{l<m}b_{k+l}e^l   &= t_{k+m}e^m.
 \end{align*}
 In each case the right hand side converges $e$--adically to $0$, so in
 $\Rhe$ we have
 \begin{align*}
  s_{ij} &= \sum_l a_{i+l,j}e^l \\
  t_k    &= \sum_l b_{k+l}e^l.
 \end{align*}
 Using this, we see that $\hat{\psi}$ is surjective.  Moreover, as a
 special case of the first equation, we have $b_j=\sum_la_{lj}e^l$.
 Putting this into the second equation gives
 $t_k=\sum_{l,m}a_{k+l,m}e^{l+m}$.  In particular, we have
 $[2](e)=\sum_{l,m}a_{lm}e^{l+m}=t_0=0$, so $\hat{\psi}$ factors
 through a map $\ov{\psi}\:\Rhat=L\psb{e}/[2](e)\xra{}\Rhe$.  As
 $\hat{\psi}$ is surjective and $\hat{\phi}\ov{\psi}=1$ we see that
 $\hat{\phi}$ is an isomorphism.
\end{proof}

\section{Inverting $e$}

\begin{proposition}\label{prop-theta}
 There is a map $\tht\:R\xra{}R'$ given by
 \begin{align*}
  \tht(s_{ij}) &= b_jb_0^{-i} - \sum_{l=1}^i a_{i-l,j}b_0^{-l} \\
  \tht(t_k)    &= -\sum_{l=1}^k b_{k-l}b_0^{-l} 
 \end{align*}
 Moreover, this induces an isomorphism $R[1/e]\xra{}R'$.
\end{proposition}
\begin{proof}
 First observe that $\tht(s_{0j})=b_j$ so the notation is
 self-consistent. 

 To see that $\tht$ respects the relations, note that
 \begin{align*}
  \tht(t_0)    &= 0 \\
  \tht(s_{10}) &= b_0b_0^{-1} - a_{00}b_0^{-1} = 1 \\
  \tht(s_{i0}) &= b_0^{1-i} -\sum_{l=1}^i a_{i-l,0} b_0^{-l} \\
               &= b_0^{i-1} - \sum_{l=1}^i \dl_{i-l,1} b_0^{-l} \\
               &= 0 \qquad \text{ for } i>1 \\
  \tht(t_k)-\tht(b_k) &= -b_k - \sum_{l=1}^k b_{k-l}b_0^{-l} \\
               &= -\sum_{l=0}^{k} b_{k-l}b_0^{-l} \\
               &= b_0(-\sum_{m=1}^{k+1} b_{k+1-m}b_0^{-m}) \\
               &= b_0\tht(t_{k+1}) 
\end{align*}  
\begin{align*}
\tht(s_{ij}) - a_{ij} &=
                  b_jb_0^{-i} - \sum_{l=0}^i a_{i-l,j} b_0^{-l} \\
               &= b_0(b_jb_0^{-1-i} - 
                  \sum_{m=1}^{i+1} a_{i+1-m,j} b_0^{-m}) \\
               &= b_0\tht(s_{i+1,j}).
 \end{align*}
 It follows that we have a ring homomorphism as described.  It induces
 a map $R[1/e]\xra{}R'$, which we again call $\tht$.  On the other
 hand, it is clear that there is a unique map $\sg\:R'\xra{}R[1/e]$
 sending $b_i$ to $b_i$, and that $\tht\sg=1$.  As
 $s_{j+1,k}=(s_{jk}-a_{jk})/e$ in $R[1/e]$, we see inductively that
 $s_{jk}$ lies in the image of $\sg$ for all $j$ and $k$.  A similar
 argument shows that $t_k$ lies in the image of $\sg$, so $\sg$ is
 surjective.  As $\tht\sg=1$ we deduce that $\tht$ and $\sg$ give an
 isomorphism $R[1/e]\simeq R'$.
\end{proof}

\begin{lemma}\label{lem-square}
 The following diagram commutes:
$$ \begin{CD}
R @>\tht>>R'\\
@V\tht VV @V\phi VV\\
\Rhat @>>\zt>(\Rhat)' 
%  \node{R} \arrow{e,t}{\tht} \arrow{s,l}{\phi}\node{R'} \arrow{s,r}{\xi} \\
%  \node{\Rhat} \arrow{e,b}{\zt} \node{(\Rhat)'}
 \end{CD}$$
\end{lemma}
\begin{proof}
 First, recall that $\xi(b_j)$ is the coefficient of $x^j$ in the
 series $e+_Fx=\sum_{m,j}a_{mj}e^mx^j$, so
 $\xi(b_j)=\sum_{m=0}^\infty a_{mj}e^m$.  Next, we have
 $\tht(s_{ij})=b_jb_0^{-i} - \sum_{l=1}^i a_{i-l,j}b_0^{-l}$.  By
 putting $m=i-l$ we can rewrite this as
 $$b_je^{-i}-\sum_{m=0}^{i-1}a_{mj}e^{m-i} $$ It follows that
 $\xi\tht(s_{ij})=\sum_{m=i}^\infty a_{mj}e^{m-i}
                 =\sum_{k=0}^\infty a_{i+k,j}e^k$.  By inspecting
 Proposition~\ref{prop-phi}, we see that this is the same as
 $\zt\phi(s_{ij})$.  

 Next, recall that $\tht(t_k)=-\sum_{l=1}^kb_{k-l}e^{-l}$, so
 \begin{align*} 
  \xi\tht(t_k) &= -\sum_{l=1}^k\sum_{m=0}^\infty a_{m,k-l}e^{m-l} \\
               &= -\sum_r e^{r-k} \sum\{a_{ij} \st i+j=r \;,\;
                                                   0\leq i<k \;,\;
                                                   0\leq j\}
 \end{align*}
 (by putting $r=m-l-k$ and $j=m$ and $i=k-l$ and noting that
 $a_{ij}=a_{ji}$).  Next observe
 \[ \sum_re^{r-k}\sum\{a_{ij}\st i+j=r \;,\; 0\leq i,j\} =
    e^{-k} [2](e) = 0 \text{ in } (\Rhat)'.
 \]
 It follows that
 \[ \xi\tht(t_k) =
     \sum_r e^{r-k} \sum\{a_{ij}\st i+j=r\;,\;i\geq k\;,\;j\geq 0\}.
 \]
 On the other hand, we have
 \[ \zt\phi(t_k)=\sum_{i,j\geq 0} a_{i+k,j}e^{i+j}. \]
 This is the same up to reindexing, as required.
\end{proof}

\section{$e$--power torsion}

Recall that $(1+t_1)e=0$ in $R$.  We shall check that this gives all
the $e$--power torsion, or equivalently that $e$ is a regular element
in the ring $\Rb=R/(1+t_1)$.  Our method is to exhibit $\Rb$ as a
colimit of rings $\Rb_k$ in which the corresponding fact can be
checked directly.

\begin{proposition}\label{prop-regular}
 The element $e$ is regular in $\Rb$.
\end{proposition}
\begin{proof}
 For $k>1$ we put
 \begin{align*}
  A_k &= L[t_k,s_{kj}\st 0<j] \\
  B_k &= L[t_i,s_{ij}\st 0\leq i\leq k,0<j].
 \end{align*}
 We define a polynomial $g_k(e)\in A_k[e]$ by
 \[ g_k(e) = \sum_{l,m=0}^{k-1} a_{ml}e^{m+l-1} +
             \sum_{l=1}^{k-1} s_{kl}e^{k+l-1} +
             t_k e^{k-1}.
 \]
 We then put $\Rb_k=A_k[e]/g_k(e)$.  One checks that $g_k(0)=2$, which
 is a regular element in $A_k$.  It follows easily that $g_k(e)$ is
 regular in $A_k[e]$, and thus that $e$ is regular in $\Rb_k$.

 Next, we define a map $\pi\:B_k\xra{}A_k[e]$ by
 \begin{align*}
  \pi(t_i) &= \sum_{l=0}^{k-1-i}\sum_{m=0}^{k-1} a_{m,i+l}e^{m+l} +
              \sum_{l=0}^{k-1-i}s_{k,i+l} e^{k+l} +
              t_k e^{k-i} \\
  \pi(s_{ij}) &= \sum_{l=0}^{k-1} a_{i+l,j}e^l + s_{kj} e^{k-i}.
 \end{align*}
 In the case $j=0$ of the second equation, the term $s_{k0}$ is to be
 interpreted as $0$.  It is easy to see that $\pi(t_k)=t_k$ and
 $\pi(s_{kj})=s_{kj}$ and $\pi(s_{00})=e$, so the notation is
 consistent.  Moreover, we see that $\pi(1+t_1)=g_k(e)$.  Using this,
 we see that $\pi$ induces an isomorphism $B_k/I_k\xra{}\Rb_k$, where 
 \begin{align*}
  I_k =& (t_0,t_1+1,s_{10}-1) + \\
       &  (s_{i0}\st i>1) + \\
       &  (t_i-s_{0i}-s_{00}t_{i+1}\st 0\leq i<k) + \\
       &  (s_{ij}-a_{ij}-s_{00}s_{i+1,j}\st 0\leq i<k,0\leq j).
 \end{align*}
 It is clear from this that $\Rb=\colim_k\Rb_k$.  As $e$ is regular in
 $\Rb_k$ for all $k$, it must be regular in $\Rb$ as well.
\end{proof}

\begin{corollary}\label{cor-bounded-torsion}
 The $e$--power torsion in $R$ is generated by $t_1+1$, and thus is
 annihilated by $e$.
\end{corollary}
\begin{proof}
 Suppose that $u\in R$ and $e^ku=0$ for some $k>0$.  It is clear from
 the proposition that the image of $u$ in $\Rb$ must be zero, so
 $u\in(t_1+1)$, so $eu=0$.
\end{proof}

\begin{corollary}\label{cor-hasse}
 There is a pullback square of rings as follows:
$$ \begin{CD}
R @>\phi>>R'\\
@V\tht VV @VVV\\
\Rhat @>>>(\Rhat)'
%  \node{R}\arrow{e,t}{\phi}\arrow{s,l}{\tht}\node{R'}\arrow{s} \\
%  \node{\Rhat}\arrow{e}\node{(\Rhat)'}
 \end{CD}$$
 Thus, $R=MU_A^*$.
\end{corollary}
\begin{proof}
 Combine Corollary~\ref{cor-bounded-torsion}, Theorem~\ref{thm-hasse},
 Proposition~\ref{prop-phi}, Proposition~\ref{prop-theta},
 Lem\-ma~\ref{lem-square}, and Theorem~\ref{thm-pullback}.
\end{proof}

\end{document}

%% file: gtoutput.tex
\def\ifplaintex{\expandafter\ifx\csname documentclass\endcsname\relax}

\ifplaintex 
\else
\fi

\expandafter\ifx\csname beginpicture\endcsname\relax
\expandafter\ifx\csname documentclass\endcsname\relax
\input pictex \else
\input prepictex \input pictex \input postpictex \fi\fi

\def\gt{{\mathsurround=0pt\it $\cal G\mskip-2mu$eometry \&\ 
$\cal T\!\!$opology}}        %  journal title in recommended style

\def\gtp{{\mathsurround=0pt\it $\cal G\mskip-2mu$eometry \&\ 
$\cal T\!\!$opology $\cal P\!$ublications}}  % GT publications

\def\lognumber#1{\def\thelognumber{#1}}
\def\volumenumber#1{\def\thevolumenumber{#1}}
\def\papernumber#1{\def\thepapernumber{#1}}
\def\volumeyear#1{\def\thevolumeyear{#1}}

\def\pagenumbers#1#2{\def\startpage{#1}\def\finishpage{#2}}
\def\published#1{\def\publishdate{#1}}
\def\proposed#1{\def\theproposer{#1}}
\def\seconded#1{\def\theseconders{#1}}
\def\received#1{\def\receiveddate{#1}}
\def\revised#1{\def\reviseddate{#1}}
\def\accepted#1{\def\accepteddate{#1}}

\def\coverauthors#1{\def\thecoverauthors{#1}}
\def\asciiauthors#1{\def\theasciiauthors{#1}}

\long\def\asciiabstract#1{\long\def\theasciiabstract{#1}}

\let\\\par\let\thelognumber\relax
\let\thevolumenumber\relax\let\thepapernumber\relax
\let\thevolumeyear\relax\let\thesamplenumber\relax\let\startpage\relax
\let\finishpage\relax\let\publishdate\relax\let\receiveddate\relax
\let\reviseddate\relax\let\accepteddate\relax\let\theasciititle\relax
\let\theasciiauthors\relax
\let\theasciiabstract\relax
\let\theasciiemail\relax\let\theshortauthors\relax\let\theshorttitle\relax
\let\thecoverauthors\relax

\long\def\maketitlep{   % start of definition of \maketitlep

\count0=\startpage

\gt\hfill      %   Journal title (top left) 
\beginpicture
\setcoordinatesystem units <0.33truein, 0.33truein> point at 2.2 0.9
\setplotsymbol ({$\cal G$})
\plotsymbolspacing=9truept
\circulararc 315 degrees from 0 1 center at 0 0
\setplotsymbol ({$\cal T$})
\circulararc 315 degrees from 1 -1 center at 1 0
\endpicture
\break
{\small\ifx\thesamplenumber\relax % sample?  
Volume \else Sample
\fi\thevolumenumber\ (\thevolumeyear)
\startpage--\finishpage\nl
Published: \publishdate}
\vglue 0.5truein plus 0.4fil minus 0.1truein

{\parskip=0pt\leftskip 0pt plus 1fil\def\\{\par\smallskip}{\ifplaintex\large
\else\Large\fi\bf\thetitle}\par\medskip}   

\vglue 0pt plus 0.1fil 

{\parskip=0pt\leftskip 0pt plus 1fil\def\\{\par}{\sc\theauthors}
\par\medskip}

\vglue 0pt plus 0.1fil 

{\small\parskip=0pt\let\newline\\
{\leftskip 0pt plus 1fil\def\\{\par}{\sl\theaddress}\par}
\expandafter\ifx\theemail\relax    % email address?
\relax\else\vglue 5pt plus 0.02fil minus 2pt\def\\{\stdspace{\rm 
and}\stdspace} 
\cl{Email:\stdspace\tt\theemail}\fi
\ifx\theurl\relax                  % URL given?
\relax\else\vglue 5pt plus 0.02fil minus 2pt\def\\{\stdspace{\rm 
and}\stdspace}
\cl{URL:\stdspace\tt\theurl}\fi\par}

\vglue 7pt plus 0.3fil minus 3pt

{\bf Abstract}
\vglue 5pt plus 0.1fil minus 2pt

\theabstract

\vglue 7pt plus 0.3fil minus 3pt

{\bf AMS Classification numbers}\quad Primary:\quad \theprimaryclass

Secondary:\quad \thesecondaryclass

\vglue 5pt plus 0.3fil minus 2pt

{\bf Keywords}\quad \thekeywords

\vglue 10pt plus 0.5fil minus 5pt

{\small  Proposed: \theproposer\hfill Received: \receiveddate\nl
Seconded: \theseconders\hfill 
\ifx\reviseddate\relax                         % paper revised?
Accepted: \accepteddate                        % no
\else
Revised: \reviseddate                          % yes
\fi}
\eject
}       %  end of definition of \maketitlep

\let\maketitlepage\maketitlep
\let\maketitle\maketitlepage

\font\phead=cmsl9 scaled 950
\font=cmsl9 scaled 1050
\font\pnum=cmbx10 scaled 913
\font\lnum=cmbx10 
\font\pfoot=cmsl9 scaled 950
\font=cmsl9 scaled 1050
\ifplaintex
\headline{\vbox to 0pt{\vskip -4.5mm\line{\small\phead\ifnum
\count0=\startpage ISSN 1364-0380 (on line)
1465-3060 (printed) \hfill {\pnum\folio}\else\ifodd\count0\def\\{ }% 
\ifx\theshorttitle\relax\thetitle\else\theshorttitle\fi\hfill{\pnum\folio}
\else\def\\{ and }{\pnum\folio}\hfill\ifx\theshortauthors\relax\theauthors
\else\theshortauthors\fi\fi\fi}\vss}}
\footline{\vbox to 0pt{\vglue 0mm\line{\small\pfoot\ifnum\count0=\startpage
\copyright\ \gtp\hfill\else
\gt, Volume \thevolumenumber\ (\thevolumeyear)\hfill\fi}\vss
}}
\else
\makeatletter
\def\@oddhead{{\small\lhead\ifnum\count0=\startpage ISSN 1364-0380 (on line)
1465-3060 (printed) \hfill {\lnum\number\count0}\else\ifodd\count0
\def\\{ }\ifx\theshorttitle\relax \thetitle \else\theshorttitle\fi\hfill
{\lnum\number\count0}\else\def\\{ and }{\lnum\number\count0}
\hfill\ifx\theshortauthors\relax 
\theauthors\else\theshortauthors\fi\fi\fi}}\def\@evenhead{\@oddhead}
\def\@oddfoot{\small\lfoot\ifnum\count0=\startpage\copyright\ \gtp\hfill\else
\gt, Volume \thevolumenumber\ (\thevolumeyear)\hfill\fi}
\def\@evenfoot{\@oddfoot}
\makeatother
\fi

\newwrite\gtoutfile
\long\gdef\makeheadfile{  %%% start of definition of \makeheadfile
{\def\\{, }\def\s{ }
\immediate\openout\gtoutfile head.xxx
\immediate\write\gtoutfile{To: math@arxiv.org}
\immediate\write\gtoutfile{Subject: put or rep NNNNN:pppp}
\immediate\write\gtoutfile{--text follows this line--}
\immediate\write\gtoutfile{Proxy-for: \ifx\theasciiauthors\relax
\theauthors\else\theasciiauthors\fi\s<\ifx\theasciiemail\relax\theemail\else\theasciiemail\fi>}
\immediate\write\gtoutfile{\noexpand\\}
\immediate\write\gtoutfile{Authors: \ifx\theasciiauthors\relax
\theauthors\else\theasciiauthors\fi}
\immediate\write\gtoutfile{Title: \ifx\theasciititle\relax
\thetitle\else\theasciititle\fi}
\immediate\write\gtoutfile{Subj-class: GT or SG or MG etc}
\immediate\write\gtoutfile{MSC-class: \theprimaryclass\ifx\thesecondaryclass\relax\else, \thesecondaryclass\fi}
\immediate\write\gtoutfile{Journal-ref: Geom. Topol. \thevolumenumber
(\thevolumeyear) \startpage-\finishpage}
\immediate\write\gtoutfile{Comments: Published by Geometry and Topology at}
\immediate\write\gtoutfile{\s\s http://www.maths.warwick.ac.uk/gt/GTVol\thevolumenumber/paper\thepapernumber.abs.html}
\immediate\write\gtoutfile{\noexpand\\}
\immediate\write\gtoutfile{}
\ifx\theasciiabstract\relax
\immediate\write\gtoutfile{\theabstract}\else
\immediate\write\gtoutfile{\theasciiabstract}\fi
\immediate\write\gtoutfile{}
\immediate\write\gtoutfile{\noexpand\\}
\immediate\write\gtoutfile{}
\immediate\closeout\gtoutfile}}  %%% end of definition of \makeheadfile

\def\maketitlepage{\maketitlep\makeheadfile}
\let\maketitle\maketitlepage

\def\ifplaintex{\expandafter\ifx\csname documentclass\endcsname\relax}

\ifplaintex 
\else
\fi

\expandafter\ifx\csname beginpicture\endcsname\relax
\expandafter\ifx\csname documentclass\endcsname\relax
\input pictex \else
\input prepictex \input pictex \input postpictex \fi\fi

\def\gt{{\mathsurround=0pt\it $\cal G\mskip-2mu$eometry \&\ 
$\cal T\!\!$opology}}        %  journal title in recommended style

\def\gtp{{\mathsurround=0pt\it $\cal G\mskip-2mu$eometry \&\ 
$\cal T\!\!$opology $\cal P\!$ublications}}  % GT publications

\def\lognumber#1{\def\thelognumber{#1}}
\def\volumenumber#1{\def\thevolumenumber{#1}}
\def\papernumber#1{\def\thepapernumber{#1}}
\def\volumeyear#1{\def\thevolumeyear{#1}}

\def\pagenumbers#1#2{\def\startpage{#1}\def\finishpage{#2}}
\def\published#1{\def\publishdate{#1}}
\def\proposed#1{\def\theproposer{#1}}
\def\seconded#1{\def\theseconders{#1}}
\def\received#1{\def\receiveddate{#1}}
\def\revised#1{\def\reviseddate{#1}}
\def\accepted#1{\def\accepteddate{#1}}

\def\coverauthors#1{\def\thecoverauthors{#1}}
\def\asciiauthors#1{\def\theasciiauthors{#1}}

\long\def\asciiabstract#1{\long\def\theasciiabstract{#1}}

\let\\\par\let\thelognumber\relax
\let\thevolumenumber\relax\let\thepapernumber\relax
\let\thevolumeyear\relax\let\thesamplenumber\relax\let\startpage\relax
\let\finishpage\relax\let\publishdate\relax\let\receiveddate\relax
\let\reviseddate\relax\let\accepteddate\relax\let\theasciititle\relax
\let\theasciiauthors\relax
\let\theasciiabstract\relax
\let\theasciiemail\relax\let\theshortauthors\relax\let\theshorttitle\relax
\let\thecoverauthors\relax

\long\def\maketitlep{   % start of definition of \maketitlep

\count0=\startpage

\gt\hfill      %   Journal title (top left) 
\beginpicture
\setcoordinatesystem units <0.33truein, 0.33truein> point at 2.2 0.9
\setplotsymbol ({$\cal G$})
\plotsymbolspacing=9truept
\circulararc 315 degrees from 0 1 center at 0 0
\setplotsymbol ({$\cal T$})
\circulararc 315 degrees from 1 -1 center at 1 0
\endpicture
\break
{\small\ifx\thesamplenumber\relax % sample?  
Volume \else Sample
\fi\thevolumenumber\ (\thevolumeyear)
\startpage--\finishpage\nl
Published: \publishdate}
\vglue 0.5truein plus 0.4fil minus 0.1truein

{\parskip=0pt\leftskip 0pt plus 1fil\def\\{\par\smallskip}{\ifplaintex\large
\else\Large\fi\bf\thetitle}\par\medskip}   

\vglue 0pt plus 0.1fil 

{\parskip=0pt\leftskip 0pt plus 1fil\def\\{\par}{\sc\theauthors}
\par\medskip}

\vglue 0pt plus 0.1fil 

{\small\parskip=0pt\let\newline\\
{\leftskip 0pt plus 1fil\def\\{\par}{\sl\theaddress}\par}
\expandafter\ifx\theemail\relax    % email address?
\relax\else\vglue 5pt plus 0.02fil minus 2pt\def\\{\stdspace{\rm 
and}\stdspace} 
\cl{Email:\stdspace\tt\theemail}\fi
\ifx\theurl\relax                  % URL given?
\relax\else\vglue 5pt plus 0.02fil minus 2pt\def\\{\stdspace{\rm 
and}\stdspace}
\cl{URL:\stdspace\tt\theurl}\fi\par}

\vglue 7pt plus 0.3fil minus 3pt

{\bf Abstract}
\vglue 5pt plus 0.1fil minus 2pt

\theabstract

\vglue 7pt plus 0.3fil minus 3pt

{\bf AMS Classification numbers}\quad Primary:\quad \theprimaryclass

Secondary:\quad \thesecondaryclass

\vglue 5pt plus 0.3fil minus 2pt

{\bf Keywords}\quad \thekeywords

\vglue 10pt plus 0.5fil minus 5pt

{\small  Proposed: \theproposer\hfill Received: \receiveddate\nl
Seconded: \theseconders\hfill 
\ifx\reviseddate\relax                         % paper revised?
Accepted: \accepteddate                        % no
\else
Revised: \reviseddate                          % yes
\fi}
\eject
}       %  end of definition of \maketitlep

\let\maketitlepage\maketitlep
\let\maketitle\maketitlepage

\font\phead=cmsl9 scaled 950
\font=cmsl9 scaled 1050
\font\pnum=cmbx10 scaled 913
\font\lnum=cmbx10 
\font\pfoot=cmsl9 scaled 950
\font=cmsl9 scaled 1050
\ifplaintex
\headline{\vbox to 0pt{\vskip -4.5mm\line{\small\phead\ifnum
\count0=\startpage ISSN 1364-0380 (on line)
1465-3060 (printed) \hfill {\pnum\folio}\else\ifodd\count0\def\\{ }% 
\ifx\theshorttitle\relax\thetitle\else\theshorttitle\fi\hfill{\pnum\folio}
\else\def\\{ and }{\pnum\folio}\hfill\ifx\theshortauthors\relax\theauthors
\else\theshortauthors\fi\fi\fi}\vss}}
\footline{\vbox to 0pt{\vglue 0mm\line{\small\pfoot\ifnum\count0=\startpage
\copyright\ \gtp\hfill\else
\gt, Volume \thevolumenumber\ (\thevolumeyear)\hfill\fi}\vss
}}
\else
\makeatletter
\def\@oddhead{{\small\lhead\ifnum\count0=\startpage ISSN 1364-0380 (on line)
1465-3060 (printed) \hfill {\lnum\number\count0}\else\ifodd\count0
\def\\{ }\ifx\theshorttitle\relax \thetitle \else\theshorttitle\fi\hfill
{\lnum\number\count0}\else\def\\{ and }{\lnum\number\count0}
\hfill\ifx\theshortauthors\relax 
\theauthors\else\theshortauthors\fi\fi\fi}}\def\@evenhead{\@oddhead}
\def\@oddfoot{\small\lfoot\ifnum\count0=\startpage\copyright\ \gtp\hfill\else
\gt, Volume \thevolumenumber\ (\thevolumeyear)\hfill\fi}
\def\@evenfoot{\@oddfoot}
\makeatother
\fi

\newwrite\gtoutfile
\long\gdef\makeheadfile{  %%% start of definition of \makeheadfile
{\def\\{, }\def\s{ }
\immediate\openout\gtoutfile head.xxx
\immediate\write\gtoutfile{To: math@arxiv.org}
\immediate\write\gtoutfile{Subject: put or rep NNNNN:pppp}
\immediate\write\gtoutfile{--text follows this line--}
\immediate\write\gtoutfile{Proxy-for: \ifx\theasciiauthors\relax
\theauthors\else\theasciiauthors\fi\s<\ifx\theasciiemail\relax\theemail\else\theasciiemail\fi>}
\immediate\write\gtoutfile{\noexpand\\}
\immediate\write\gtoutfile{Authors: \ifx\theasciiauthors\relax
\theauthors\else\theasciiauthors\fi}
\immediate\write\gtoutfile{Title: \ifx\theasciititle\relax
\thetitle\else\theasciititle\fi}
\immediate\write\gtoutfile{Subj-class: GT or SG or MG etc}
\immediate\write\gtoutfile{MSC-class: \theprimaryclass\ifx\thesecondaryclass\relax\else, \thesecondaryclass\fi}
\immediate\write\gtoutfile{Journal-ref: Geom. Topol. \thevolumenumber
(\thevolumeyear) \startpage-\finishpage}
\immediate\write\gtoutfile{Comments: Published by Geometry and Topology at}
\immediate\write\gtoutfile{\s\s http://www.maths.warwick.ac.uk/gt/GTVol\thevolumenumber/paper\thepapernumber.abs.html}
\immediate\write\gtoutfile{\noexpand\\}
\immediate\write\gtoutfile{}
\ifx\theasciiabstract\relax
\immediate\write\gtoutfile{\theabstract}\else
\immediate\write\gtoutfile{\theasciiabstract}\fi
\immediate\write\gtoutfile{}
\immediate\write\gtoutfile{\noexpand\\}
\immediate\write\gtoutfile{}
\immediate\closeout\gtoutfile}}  %%% end of definition of \makeheadfile

\def\maketitlepage{\maketitlep\makeheadfile}
\let\maketitle\maketitlepage